\documentclass[11pt]{article}
\usepackage{amsmath,amssymb,amsthm, mathrsfs}

\begin{document}

\renewcommand{\baselinestretch}{1.3}
\renewcommand{\arraystretch}{1.3}

\begin{center}{\bf\LARGE  Five Theorems on Splitting Subspaces and Projections in Banach Spaces and Applications to  Topology and Analysis  in Operators }\\
\vskip 0.5cm Ma Jipu$^{1,2}$

1 Department of Mathematics, Nanjing University, Nanjing 210093,
China\\

2 Tseng Yuanrong Functional Research Center, Harbin Normal University,\\ Harbin 150080, China \\

E-mail: jipuma@126.com
\end{center}

{\bf Abstract}\quad  Let $B(E,F)$ denote the set of all bounded linear operators from  $E$
into $F$, and $B^+(E,F)$ the set of  double splitting operators in
$B(E,F)$. When both $E,F$ are infinite dimensional , in $B(E,F)$ there are not more elementary transformations in matrices so that lose the way to discuss the path connectedness of such sets in $B^+(E,F)$ as $\Phi_{m,n}=\{T\in B(E,F): \dim N(T)=m<\infty \ \mbox{and} \  \mathrm{codim}R(T)=n<\infty\},$ $F_k=\{T\in B(E,F): \mathrm{rank}\, T =k<\infty\}$, and so forth. In this paper we present five theorems on projections and splitting subspaces in Banach spaces instead of the elementary  transformation. Let $\Phi $ denote any one of $F_k ,k<\infty$ and $\Phi_{m,n}$ with either $m>0$ or $n>0.$ Using these theorems we prove $\Phi$ is path connected.Also these theorems bear an equivalent relation in $B^+(E,F)$, so that the following general result follows: the equivalent class $\widetilde{T}$ generated by $T\in B^+(E,F)$ with either $\dim N(T)>0$ or $\mathrm{codim} R(T)>0$ is path connected. (This equivalent relation in operator topology appears  for the first time.) As applications of the theorems we give that $\Phi $ is a smooth and path connected submanifold in $B(E,F)$ with the tangent space $T_X\Phi =\{T\in B(E,F): TN(X)\subset R(X)\}$  at any $ X\in {\Phi },$ and prove that $B(\mathbf{R}^m,\mathbf{R}^n)=\bigcup^{\min\{n,m\}}\limits_{k=0}F_k $ possesses the following properties of geometric and topology :
$F_k ( k <\min\{ m,n\})$ is a smooth and path connected subhypersurface in $B(E,F)$, and specially, $\dim F_k=(m+n-k)k, k=0,1, \cdots , \min\{m.n\}.$ Of special interest is the dimensional formula of $F_k \, \, k=0,1, \cdots , \min\{m.n\},$ which is a new result in algebraic geometry.
In view of the proofs of the above theorems it can not be too much to say that Theorems $1.1-1.5$ provide the rules of finding path connected sets in $B^+(E,F).$



{\bf 1\quad Five Theorems on Splitting Subspaces and Projections in Banach Spaces  }

 In this section we introduce  five theorems on splitting subspaces and projections in Banach spaces instead of  the elementary  transformation in  matrices.
Let  $E,F$ be two Banach spaces, $B(E,F)$  the set of all linear
bounded operators from  $E$ into $F$ , and $B^+(E,F)$  the set of all
double splitting operators in $B(E,F)$ . Recall that two operators $T_0$ and $T_1$ in a set $S\subset B(E,F)$ are said to be path connected in $S$ provided there is a path $\omega (t)\in S$ for $t\in [0,1]$ such that $\omega (0)=T_0$ and $\omega (1)=t_0 $ ; a set $S$ in $B(E,F)$ is said to be path connected provided arbitrary two operators in $S$ are path connected in it. For simplicity, write two operators $T_0$ and $T_1$ being path path connected in $S$ as $T_0\thicksim T_1$ in $S$ .

{\bf Theorem 1.1}\quad{\it  If $E = E_{1}\oplus R = E_{*}\oplus R$ , then the following conclusions hold:

$(i)$ there exists a unique $\alpha\in B(E_*, R)$ such that
$$E_1 = \left\{x + \alpha x: \forall x\in E_*\right\};\eqno(1.1)$$
conversely, for any $\alpha\in B(E_*, R)$  the subspace $E_1$ defined by $(1.1)$ satisfies $E = E_1\oplus R $

$(ii)$ $$ P^R _{E_1} = P _{ E_*}^R + \alpha P^R _{E_*} \,\,\mbox{and  so} \,\,P_R^{E_1} = P_R ^{E_*} - \alpha P^R _{E_*} .$$ }

{\bf Proof}\quad  For the proof of $(i)$ refer [Ma3] and [Abr].
Obviously,
$\left(P^{R}_{E_*} + \alpha P^{R}_{E_*}\right)^{2} = P^{R}_{E_*} + \alpha P^{R}_{E_*},$ $ R( P _{ E_*}^R + \alpha P^R _{E_*})=E_1$ by (1.1) \, \mbox{and}
$$P^{R}_{E_*}x + \alpha P^{R}_{E_{*}}x = 0\,\Leftrightarrow P^{R}_{E_{*}}x=0
\Leftrightarrow x\in R , \, i.e., N( P _{ E_*}^R + \alpha P^R _{E_*})=R .$$
 This infers
$$
P^{R}_{E_1} = P^{R}_{E_*} + \alpha P^{R}_{E_*} ,\,\,\mbox{and so ,} \,\,P_{R}^{E_1} = P_{R}^{E_*} - \alpha P^{R}_{E_*} .
$$
The proof ends . $\Box $

{\bf Theorem 1.2}\quad{\it  Suppose that $E = E_{*}\oplus R $ and
$\dim R > 0  .$  Then $P^{R}_{E_{*}}$ and $(-P^{R}_{E_{*}})$ are
path connected in the set $ S=\{T\in B(E): E=R(T)\oplus R$ and
$N(T)=R\}$.}

{\bf Proof}\quad  Due to $\mbox{dim} R > 0,$ it is clear that $B(E_{*},R)$ contains a non-zero operator $\alpha$ if \ $\mbox{dim} E_* >0,$ otherwise
$E_*= \{0\}$ and so $P^{R}_{E_{*}}=o ,$ the theorem is trivial . We assume both , $\dim E ,\, \dim R >0 $ in the sequel . Let $E_{1} = \left\{x+ \alpha
x:\forall x\in E_{*} \right\}.$ Then by Theorem 1.1 , $E =
E_{1}\oplus R $ , and
$$
P^{R}_{E_{1}} = P^{R}_{E_{*}} + \alpha P^{R}_{E_{*}} . \eqno(1.2)
$$
Consider the path
$$\omega(t)= (1-2t)
P^{R}_{E_{*}}+(1-t)\alpha P^{R}_{E_{*}}\,\,\,\, t\in [0,1] .$$
Clearly , $\omega (0)=P^R_{E_1},\omega(1)=-P^R_{E_*}$ and
$$R(\omega(t))= R(P^{R}_{E_{*}} + \frac{1-t}{1-2t}\alpha
P^{R}_{E_{*}}) ,\quad\,\,\,t\in [0,1] .$$
Then by Theorem 1.1, $R(\omega(t))\oplus R = E$ for all $t\in
[0,1]$ . In addition , $\omega(0)=P^R_{E_1} , \, \omega(1)= -P^R_{E_*}$ and
 $ N(\omega(t))=R \, , \, \forall t\in[0,1]$ . This shows that
$P^{R}_{E_{1}}$ and $-P^{R}_{E_{*}}$ are path connected in $S$ . Next go to show that $P^{R}_{E_{1}}$ and $P^{R}_{E_{*}}$ are path
connected in $S .$  Consider the path
$$\omega(t) =  P^{R}_{E_{*}} + t\alpha P^{R}_{E_{*}}\,\,\,\, \forall t\in [0,1] .$$
By Theorem 1.1 , $\omega(1) =  P^{R}_{E_{*}} + \alpha P^{R}_{E_{*}}=  P^{R}_{E_{1}} ,$
 and $R(\omega(t))\oplus R = E \,\,\,\,\forall \in [0,1] .$
Obviously , $N(\omega(t))= R $ and $\omega(0)= P^{R}_{E_{*}} .$  So
$P^{R}_{E_{1}}$ and $P^{R}_{E_{*}}$ are path connected in
$S.$  The theorem is proved . \quad $\Box$

{\bf Theorem 1.3}\quad{\it Suppose that  $T_{0}\in B^+(E,F)$ with $ R(T_0)\oplus N=F ,$
 $F=F_*\oplus N$ and $\mbox{dim}N>0$. Then $T_{0}\sim P_{F_*}^N T_0 $
 in the set  $S=\{T\in B^+(E,F):F=R(T)\oplus N \ \mbox{ and} \  N(T)=N(T_0)\}.$}

 {\bf Proof}\quad  One can assume $R(T_{0}) \neq F_{*}$ , otherwise $P_{F_*}^N T_0= T_0 $ the
theorem is trivial. Thus by Theorem 1.1 , there exists a non-zero
operator $\alpha\in B(F_{*},N)$ such that
$$R(T_{0})=\left\{y+\alpha y : \forall y\in F_{*}\right\},\,\,\,P^{N}_{R(T_{0})}= P^{N}_{F_{*}}+\alpha P^{N}_{F_{*}} ,$$
and
$P_{N}^{R(T_{0})}= P_{N}^{F_{*}}-\alpha P^{N}_{F_{*}} .$ So
$$T_{0}=\left\{P^{N}_{F_{*}}+\alpha P^{N}_{F_{*}}\right\}T_{0} .\eqno(1.3)$$
Let
$$F_t=\left\{y+ t\alpha y : \forall y\in
F_{*}\right\}\,\,\,\,\mbox{for all} \, \, \, t\in [0,1] .$$
Due to $t\alpha \in B(F_*,N)$ for any $t\in [0,1] ,$ Theorem 1.1 shows $F_t\oplus N=F, P_{F_t}^N=P_{F_*}^N +t \alpha P_{F_*}^N $ and $P_N^{F_t}=P_N^{F_*}-t\alpha P_{F_*}^N \forall t\in [0,1]$ .\\ Consider the path $$ \omega (t)= P_{F_t}^N T_0 \quad \quad  \forall t\in [0,1] .$$ We claim the path $\omega (t)$ lies in $S$ and connects $T_0$ with $P_{F_*}^NT_0$ . Clearly$$ x\in N(\omega (t))\Leftrightarrow  P_{F_t}^N T_0 x =0\Leftrightarrow T_0x\in N \Leftrightarrow x\in N(T_0)$$
because of $R(T_0)\oplus N=F$ ( which implies $R(T_0)\cap N=\{0\}$ ), and $$R(\omega(t))=P_{F_t}^N R(T_0) =P_{F_t}^N ( R(T_0) \oplus N)=P_{F_t}^N F=F_t.$$ So $R(\omega (t))\oplus N(=F_t\oplus N)=F .$ Meanwhile $\omega(1)=T_0$ by $(1.3) ,$  \mbox{and}  $\omega (0)=P_{F_*}^N T_0.$ The proof ends .\quad $\Box$

{\bf Theorem 1.4}\quad {\it Suppose that $T_{0}\in B^+(E,F)$ with $N(T_0)\oplus R_0=E , E=E_{*}\oplus R_{0} \, \mbox{and}\, \dim R_0>0 .$ Then $T_{0}\sim TP_{R_0}^{E_*}$ in the set  $S=\{T\in B^+(E,F): N(T)\oplus R_0=E  \  \mbox{and}  \,  R(T)=R(T_0)\} .$}

\textbf{Proof}\quad  One can assume $E_{*}\neq N(T_{0})$ . Otherwise $T_0=TP_{R_0}^{E_*} ,$
the theorem is trivial.  By Theorem 1.1, there exists a non-zero
operator $\alpha\in B(E_{*},R_0)$ such that
$$P_{R_0}^{N(T_0)}=P_{R_0}^{E_*}-\alpha P^{R_0}_{E_*},$$ and so  $$ T_0=T_0P_{R_0}^{N(T_0)}=T_0(P_{R_0}^{E_*}-\alpha P_{E_*}^{R_0}). \eqno (1.4) $$
Consider the path $$\omega (t)=T_0(P_{R_0}^{E_*}-t \alpha P_{E_*}^{R_0} ) \, \,\, \forall t\in [0,1].$$ We claim $R(\omega (t))=R_0$ and $N(\omega (t)) \oplus R_0=E$. Clearly  $\omega (t)R_0= T_0(P_{R_0}^{E_*}-t \alpha P_{E_*}^{R_0})R_0=T_0R_0=R(T_0).$ Next go to verify $N(\omega (t))\oplus R_0=E.$ Let $$N_t=\{x+ t\alpha x: \forall x\in E_*\} \, \, \forall t\in [0,1].$$ By Theorem $1.1$
$P^{N_t}_{R_0}=P_{R_0}^{E_*}- t\alpha P_{E_*}^{R_0} $ Obviously
\begin {eqnarray*} x\in N(\omega(t))
&\Longleftrightarrow & T_0P^{N_t}_{R_0}x= 0 \Leftrightarrow P_{R_0}^{N_t} x\in N(T_0)\\ &\Leftrightarrow & P_{R_{0}}^{N_t}x=0 \,\, \mbox{because \,of } N(T_0)\oplus R_0=E \\ &\Leftrightarrow & x\in N_t,\end{eqnarray*} and so $N(\omega (t))\oplus R_0 ( =N_t\oplus R_0)=E  \, \mbox{for all }  \, t\in [0,1].$ By $(1.4)$ it is direct that $\omega (0)=T_0P_{R_0}^{E_*}$ and $\omega(1)=T_0.$ Then the proof of the theorem ends .\quad $\Box$

{\bf Theorem 1.5}\quad{\it Suppose that two subspaces $E_1$ and
$E_2$ in $E$ satisfy $\dim E_1=\dim E_2< \infty.$ Then $E_1$
and $E_2$  possess a common complement R, i.e.,$E=E_1\oplus R=
E_2\oplus R $.}

\textbf{Proof}\quad  Because $dim E_1=\dim
E_2<\infty$, there exit subspaces $H_1, E_1^* \, \mbox{and} \, E_2^*$ such that
$$E=H_1\oplus (E_1+E_2),\,E_1=E^*_1\oplus(E_1\cap E_2),\, \,\mbox{and} \, \  E_2=E^*_2\oplus(E_1\cap E_2),$$ respectively.
It is easy to observe that  $(E^*_1\oplus
E^*_2)\cap(E_1\cap E_2)=\{0\}$ and $\mbox{dim} E^*_1= \mbox{dim}
E^{*}_{2}<\infty.$ Indeed, if $e^*_1+e^*_2$ for $e^*_1\in E^*_1 \, \mbox{and} \, e^*_2\in E^*_2$ belongs to
$E_1\cap E_2$, then
$e_2^*=(e^*_1+e^*_2)-e^*_1\in E_1$ and
$e_1^*=(e^*_1+e^*_2)-e^*_2\in E_2$, and so
$e^*_1=e^*_2=0$ because of $E^*_i\cap(E_1\cap
E_2)=\{0\},\,i=1,2.$ Hereby,$$E_1+ E_2=(E^*_1\oplus E^*_2)\oplus(E_1\cap
E_2).\eqno(1.5)$$
We are now in the position to determine $R$. We may assume $dim
E^*_1=dim E^*_2>0,$ otherwise $E_1=E_2$, the theorem is trivial. Thus
$B^\times(E^*_1,E^*_2)$ is
 non empty. Let $\alpha \in B^\times(E^*_1,E^*_2) $ and
$$H=\{x+\alpha x:\forall x\in E^*_1\}=\{x+\alpha^{-1}x:\forall x\in E^*_2\}.$$
By Theorem 1.1,
$$E^*_1\oplus E^*_2=H\oplus E^*_ 1=H\oplus E^*_2.$$
By (1.5)
$$
E=H_1\oplus(E_1+ E_2)=H_1\oplus(E^*_1\oplus E^*_2)\oplus(E_1\cap E_2)=H\oplus H_1\oplus E_2^*\oplus(E_1\cap E_2)=H\oplus H_1\oplus E_2
.$$
Similarly $$
E=H\oplus H_1\oplus E_1^*\oplus(E_1\cap E_2)=H\oplus H_1\oplus E_1.$$
This infers  $R=H\oplus H_1.$ The proof ends.\quad$\Box$

\begin{center}{\bf 2\quad Some Path Connected Sets in $B ^+(E,F)$}
\end{center}\vskip 0.2cm

 In this section we apply The five theorems to study  path connectedness of such sets in $B^{\times}(E,F)$ as $F_k, \, \Phi_{m,n}$ and forth,
where $F_k =\{ T\in B(E,F): rankT=k< \dim  F \}$ and $\Phi_{m,n}=\{T\in B^+(E,F):\dim N(T)=m \, \, \mbox{and} \, \, \mathrm{codim}\,R(T)=n\}.$
Let $\Phi$ denotes any one of $ F_k  \, \, \mbox{and} \, \, \Phi _{m,n} $  with either  $  n >0 $  or  $  m >0.$

{ \bf Theorem $2.1$ }\quad{\it The following conclusion for any $T\in \Phi $ holds: $T\sim -T$ in $\Phi $ .}

{\bf Proof}\quad  We first consider  $\Phi = F_k  \,\, \mbox{or} \,\, \Phi_{m,n} $ with $ n >0 ,$ say $$R(T)\oplus N^+=F  \, \, \mbox{and} \, \, N(T)\oplus R=E \, \, \, \mbox{for } \, \, T\in \Phi.$$  By assumption  $ \dim N^+>0.$  Then by Theorem $1.2$ there exist a path $\omega (t)$ lying in the set $S=\{ T\in B(F):N(T)=N^+  \, \, \mbox{and}  \, \, R(T)\oplus N^+=F \} $ i.e., $N(\omega (t))=N^+ \, \mbox{and} \, R(\omega (t))\oplus N^+=F,$ such that $\omega (0)=P_{R(T)}^{ N^+ } \, \, \mbox{and} \, \, \omega (1)=-P_{R(T)}^{N^+ }.$ \\
 Consider the path as follows, $\gamma(t)=\omega(t)T.$ We claim that $\gamma(t)$ is just such a path lying in $\Phi$ that $T\sim -T$ in $\Phi.$ Indeed,$$R(\gamma(t))=\omega(t)R(T)=\omega(t)(R(T)\oplus N^+)=\omega(t)F=R(\omega(t)) $$ and  $$ x\in N(\gamma (t))\Leftrightarrow Tx\in N(\omega (t))\Leftrightarrow Tx \in N^+ \Leftrightarrow x\in N(T) \mbox{because of } \, R(T)\oplus N^+=F . $$
 So $$ R(\gamma(t))\oplus N^+= R(T)\oplus N^+=F ,\eqno (2.1)$$ and  $$ N(\gamma(t))\oplus R= N(T)\oplus R=E.\eqno (2.2) $$ Then the following conclusion for $T\in F_k$  follows from $(2.1)$ :  $\gamma (t)$  for $t\in [0,1]$ lies in $F_k$ such that $\gamma(0)=T \, \mbox{and} \, \gamma(1)=-T ;$ similarly, since $(2.1) \, \, \mbox{and} \, \, (2.2)  $ we infer for $T\in \Phi_{m,n} $ with $ n>0,$  $\gamma (t)\in \Phi_{m,n} $ with $n>0 \, \, \forall t\in [0,1]$. So the theorem for $\Phi=F_k  \, \, \mbox{or} \, \,  \Phi_{m,n} $ with $ n>0$  holds.  The residue is to prove that the theorem for $\Phi_{m,0}$ with $m>0$ holds. Let $E=N(T)\oplus R^+$ for $T\in \Phi_{m,0} $ with $m > 0.$  Due to $\dim N(T)=m>0,$ by Theorem $1.2$ there exists a path $\omega (t)$ with $\omega (0)=P_{R^+}^{N(T)} $ and  $ \omega (1)=-P_{R^+}^{N(T)}$  satisfying  that $N(\omega (t))=N(T)$ and $R(\omega (t))\oplus N(T)=E . $ Write $\gamma (t)=T\omega (t).$  Obviously
$$R(\gamma(t))=TR(\omega (t))=T(R(\omega (t)) \oplus N(T))=TE=F $$ and
$$\gamma(t)x=0\Leftrightarrow \omega(t)x\in N(T)\Leftrightarrow x\in N(T) \, \, \mbox{because} \, \, R(\omega (t))\oplus N(T)=E.$$ So $\gamma (t)$ lies in $\Phi_{m,0}$ with $m>0$ such that $T \sim -T$ in $\Phi_{m,0}$.The proof ends.\quad $\Box$

{\bf Remark } \quad  The  theorem for $\Phi_{0,0}$ is false. For example $\Phi_{0,0}=B^\times (F).$  Theorem $2.1$ is very interesting. It will play a crucial role in the next proofs of Theorems $2.2 - 2.4$.

{\bf Theorem $2.2$ } \quad {\it $\Phi_{m,n}$ for $m,n<\infty$ and either $n>0$ or $m>0$  is path connected.}

{\bf Proof } \quad First go to show that $\Phi_{m,0}$ with $\infty >m
>0$ is path connected. Let $T_1 \, \mbox{and} \, T_2$ be arbitrary two operators in  $\Phi_{m,0}.$
 Because $\dim N(T_1)=\dim N(T_2)= m <\infty,$ by Theorem $1.5$ there exists a subspace $R$ in $E$ such that the following conclusions for $T_1$ and $T_2\in \Phi_{m,0}$  hold:
$$ E=N(T_i)\oplus R \quad \mbox{and} \quad T_i|_R \in B ^{\times} (R,F) \quad i=1, 2 .\quad \eqno (2.3)$$
Let $T_1^+= (T_1|_ R )^{-1}\in B(F,R).$ By Theorem $1.4$ it is easy to see  $$ T_2\sim T_2P_R^{N(T_1)}=T_2T_1^+T_1 \quad  \mbox{in} \, \, \Phi_{m,0} \, \, \mbox{with} \, \, m >0.$$
Meanwhile $T_2T_1^+\in B^\times (F).$  Indeed,  $R(T_2T_1^+)=T_2R(T_1^+)=T_2R=F$ and
$$ y\in N(T_2T_1^+)\Leftrightarrow T_1^+y\in N(T_2)\Leftrightarrow T_1^+y=0 \,\, \mbox{by $(2.3)$} \Leftrightarrow y=0 \,\, \mbox{by $(2.3)$}.$$
Thus $T_2T_1^+\sim I_F $ or $-I_F$ in $\Phi_{m,0},$ so that $$T_2\sim T_2T_1^+T_1 \sim T_1 \, \, \, \mbox{or} \, \, \,- T_1 \quad \mbox{in}  \quad \Phi_{m,0}.$$
By Theorem $2.1$ we conclude  $T_2 \sim T_1$ in $\Phi_{m,0}.$\\
 Finally to prove that $\Phi_{m,n}$ with $n>0$ is path connected. By Theorem $1.5$ there exists a subspace $F_0$ such that the following equalities for $T_1$ and $T_2 \in \Phi_{m,n}$ with $n>0$ hold: $$F=F_0 \oplus N_1=R(T_1) \oplus N_1 \quad \mbox{and} \quad F=F_0\oplus N_2= R(T_2) \oplus N_2$$ because $\dim N_1=\dim N_2<\infty.$ Then by Theorem $1.3$ $$P_{F_0}^{N_i}T_i\sim T_i \quad \mbox {in} \quad \Phi_{m,n} \quad i=1,2 .$$ For simplicity write $P_{F_0}^{N_i}T_i$ as $T_i^0 \quad i=1,2.$ We claim that $R(T_i^0)=F_0$ and $N(T_i^0)=N(T_i) \quad i=1,2 .$ In fact $$R(T_i^0)=P_{F_0}^{N_i}R(T_i)=P_{F_0}^{N_i}(R(T_i)\oplus N_i)=P_{F_0}^{N_i}F=F_0$$ and $$T_i^0x=0\Leftrightarrow T_ix\in N_i\Leftrightarrow T_ix=0 \, \, \mbox{because} \, \, R(T_i)\cap N_i=\{0\}\Leftrightarrow x\in N(T_i) \quad i=1,2.$$ Thus the proof of the theorem turns to that of $T_i^0\sim T_2^0$ in $\Phi_{m,n}.$
Let  $\Phi_{m,0}(E,F_0)=\{T\in B(E,F_0):\dim N(T)=m \, \, \mbox {and} \, \,\mathrm{codim} R(T)=0\}$. Obviously $\Phi_{m,0}(E,F_0)\subset \Phi_{m,n},$  and so , $T_1^0 \sim T_2^0$ in $\Phi_{m,n}$ if $T_1^0\sim T_2^0$ in $\Phi _{m,0}(E,F_0).$ Take the place of $F$ in the theorem by $F_0$, then by the proved theorem for $\Phi_ {m,0}$ ,  $T_1^0\sim T_2^0$ in $\Phi_{m,0}(E,F_0)$ and so, $T_1^0\sim T_2^0$ in $\Phi_{m,n}.$
The theorem is proved .$ \quad \Box $

{\bf Theorem $2.3$}\quad {\it $F_k \, \mbox{for} \, k< \dim F $ is path connected .}

{\bf Proof }\quad  The theorem for $k=0$ is trivial. In what follows, we assume $k>0$.
Let $T_1, T_2$ be arbitrary two operators in $F_k,$ say $E=N(T_i)\oplus R_i, \, \, \mbox{and} \, \,  F=R(T_i)\oplus N_i^+, i=1, 2$. Obviously $\dim R_i=
 \dim R(T_i)= k<\infty \, \, i=1,2 .$ By Theorem $1.5$ there are exist subspace $N_0 $ in $E$ and a subspace $N^+$ in $F$ such that $$E=N_0\oplus R_1=N_0\oplus R_2\, \mbox{and} \,  F=R(T_1)\oplus N^+ =R(T_2)\oplus N^+ .$$
Let
$$L_ix=\left\{\begin{array}{rcl}T_ix, & & {x\in R_i}\\
0, & & {x\in N_0}
\end{array}\right.$$
$i=1,2.$
We claim $T_i\sim L_i$ in $F_k \,\, i=1,2$. Since $E=N(T_i)\oplus R_i = N_0\oplus R_i \, \, i=1, 2 ,$ by Theorem $1.1$  there exists an operator $\alpha _i\in B(N_0,R_i) $ such that
$$P^{N(T_i)}_{R_i}=P_{R_i}^{N_0} - \alpha _i P_{N_0}^{R_i} \, \, i=1,2. \quad \eqno(2.4)$$ Consider the following path $\omega _i(t) \, \,  i=1,2 :$
$$\omega_i(t)=T_i(P_{R_i}^{N_0} - t\alpha _iP_{N_0}^{R_i}) \, \,\, \forall \, \, t\in [0,1].$$ Such path $\omega_i(t)$ lies in $F_k$ that $T_i\sim L_i$ in $F_k$. Indeed ,  $\omega_i(0)=T_iP_{R_i}^{N_0}=T_i$ and $\omega_i(1)=T_iP_{R_i}^{N(T_i)}=L_i$ by $(2.4)$ ; while from $R(\omega_i(t))\subset T_i R_i \, \mbox{and} \,  \omega_i(t)x= T_i x  \, \, \forall \, \, x\in R_i$ it follows that $R(\omega_i(t))= T_iR_i$ for $t\in [0,1] , i=1,2$. So $T_i\sim L_i$ in $ F_k.$ Thus the proof of the theorem turns to that of $L_2\sim L_1$ in $F_k.$ \\
Note that $R(L_1)(=R(T_1)) \oplus N^+ = R(L_2)(=R(T_2))\oplus N^+=F \, \mbox{and} \, \dim N^+>0 .$ It is immediate $L_2\sim P_{R(L_1)}^{N^+} L_2$ in $F_k$ from Theorem $ 1.3 .$ Hence, it is enough to show $P_{R(T_1)}^{N^+}L_2\sim L_1 \, \mbox{in} \, F_k .$ Let $$L_1^+y=\left\{\begin{array}{rcl}(L_1|_{R_1})^{-1}y, & & {y\in L_1}\\
0, & & {y\in N^+}.
\end{array}\right .$$
 Then $$P_{R(L_1)}^{N^+}L_2=P_{R(L_1)}^{N^+}L_2P_{R_1^+}^{N_0}=P_{R(L_1)}^{N^+}L_2L_1^+ L_1 \quad \eqno (2.5)$$ (note $N(L_2)=N_0$ and $P_{R_1^+}^{N_0}=L_1^+ L_1).$ Especially, $$P_{R(L_1)}^{M^+}L_2L_1^+|_{R(L_1)}\in B^\times (R(L_1)).\quad \eqno (2.6)$$ In fact \begin {eqnarray*} && P_{R(L_1)}^{N^+}L_2L_1 ^+y=0 \, \,  \mbox{for} \, \, y\in R(L_1) \Leftrightarrow L_2L_1^+y\in N^+\\ && \Leftrightarrow L_1^+y\in N_0 \, \, \mbox{since }\, \, N^+\cap R(L_2)=\{0\} \\ && \Leftrightarrow y=0 \, \, \mbox{since} \, \, R(L_1) (=R(T_1))\cap N_0=\{0\};\end{eqnarray*} while we have for any $y\in R(L_1)$
\begin{eqnarray*}  y &=&P_{R(L_1)}^{N^+}y=P_{R(L_1)}^{N^+}(P_{R(L_2)}^{N^+}y +P_{N^+}^{R(L_2)}y)\\&=&P_{R(L_1)}^{N^+}P_{R(L_2)}^{N^+}y=P_{R(L_1)}^{N^+}x \, \,( x= P_{R(L_2)}^{N^+}y \in R_2^+ )\\&=&P_{R(L_1)}^{N^+}L_2
(P_{R_1^+}^{N_0}x +P_{N_0}^{R_1^+}x)=P_{R(L_1)}^{N^+}L_2P_{R(L_1)}^{N^+}\\&=&P_{R(L_1)}^{N^+}L_2L_1^+y_1 \, \, \mbox{for\,some} \, \, y_1\in R(L_1).\end{eqnarray*} This shows that the conclusion $(2.6)$ holds. Therefore $$P_{R(L_1)}^{N^+}L_2L_1^+|_{R(L_1)}\sim I_{R(L_1)} \, \mbox{or} \,-I_{R(L_1)} \, \mbox{in} \, \,B^\times(R(l_1)),$$ Where $I_{R(L_1)}$ denotes the identity on $R(L_1).$ By $(2.5 )$ we infer $$P_{R(L_1)}^{N^+}L_2 =P_{R(L_1)}^{N^+}L_2L_1^+L_1\sim L_1 \, \mbox{or} \, -L_1 \,  \mbox{in} \, \, F_k.$$ By Theorem $2.1$  one can conclude in either way
 $$P_{R(L_1)}^{N^+}L_2L_1^+L_1  \,  i.e.,  P_{R(L_1)}^{N^+}L_2 \sim L_1 \,\,\mbox{in} \, \, \,F_k.$$ The proof ends .\quad $\Box$

let $G(.)$ denote the set of all splitting subspaces in the Banach space in the parentheses, and $$ U_E(R)=\{H\in G(E): E=H\oplus R \} \, \,\mbox{and} \, \, U_F(S)=\{L\in G(F): F=S\oplus L\}.$$
 The proofs of Theorems $2.2 \,  \mbox{and} \, 2.3$ show that the theorems are able to be proved with Theorems $1.1 - 1.4$
due to Theorem $1.5$. But Theorem $1.5$ for two subspaces of infinity dimension is false. Nevertheless , instead of it , we introduce
 the equivalent relation as follows,

{\bf Definition $2.1$} \quad   $T_0$ and $ T_*$ in $ B^+ (E,F)$ are said to be equivalent provided there exist finite number of $N_i $ in $G(E)\, \, , i=1,2, . . .,m,$ and $F_j$ in $G(F) \, \, , j=1,2, . . . ,n$ such that all
$$U_E (N(T_0))\cap U_E(N_1), \, . \, . \, .\, ,U_E(N_m)\cap U_E(N(T_*)) \quad \eqno (2.7)$$ and $$U_F(R(T_0))\cap U_F(F_1), \, . \, . \,  . \, , U_F(F_m)\cap U_F(R(T_*))\quad \eqno(2.8) $$ are non empty.\\
(This equivalent relation in operator topology appears for the first time .)\\
Let $ \widetilde{T}$ denote the equivalent class generated by $T$ in $B^+ (E,F).$ We have

 {\bf Theorem $2.4$}\quad {\it $ \widetilde{T_0} $ for any $T_0\in B^+ (E,F)$ with either $\dim N(T_0)>0$ or $\mathrm{codim}\,R(T_0)>0$ is path connected.}

{\bf Proof}\quad First we discuss the case of  $  \mathrm{codim}\,R(T_0) > 0.$ Assume that $T_{*}$ is any operator in $\widetilde{T_{0}}$, and
$$ R_1\in U_E(N(T_0))\cap U_E(N_1),
\cdots,  R_{m+1}\in U_E(N_m)\cap U_E( N(T_*)).\eqno(2.9)$$
We define inductively
$$ T_k = T_{k-1}P^{N_k}_{R_k},\,\,k=1,2,\cdots,m.$$
Evidently
$$ R(T_k) = R(T_0) \, \, \mathrm{and } \, \, N(T_k) = N_k, k=1,\cdots,m.\eqno(2.10)$$
( Indeed by $(2.9),$  $R(T_k)=T_{k-1}R_k = T_{k-1} (R_k \oplus N(T_k)( = N_k)) =T_{k-1} E=R(T_{k-1}) \, \mbox{and} $
$T_k x=0\Leftrightarrow P_{R_k}^{N_k}x\in N(T_{k-1})\Leftrightarrow x \in N_k \, \mbox{because of} \, \, E=N(T_{k-1})\oplus R_k \, \mbox{for} \, \, k=1, \cdots, m .$ Hereby $(2.10)$ follows.)

Due to $E=R_k\oplus N_{k-1}=R_k\oplus N_k, i.e., R_k\in U_E(N_{k-1})\cap U_E(N_k),$
Theorem $1.4$ shows $T_{k-1} \sim T_{k-1} P^{N_k}_{R_k}, \, \, i.e., T_k \sim T_{k-1},$ in the set  $\Omega_k= \{T\in B^+(E,E): E=N(T)\oplus R_k \, \mbox{and} \, R(T)=R(T_{k-1})\} \, \mbox{for} \, k=1, \cdots, m ;$  while by $(2.9) \, \mbox{and} \,  (2.10)$, $R_1\in U_E(N(T_0))\cap U_E(N_1) , \cdots, R_k\in U_E(N_k)\cap U_E(N(T))$. Hence $\Omega_k \subset \widetilde{T_0} \, \mbox{and}$ so $$T_k \sim T_0  \, \,\mbox{in} \, \, \widetilde{T_0} \quad \mbox{for} \, \, \,  k=1, \cdots, m . \eqno (2.11)$$
 Meanwhile assume
$$ S_1\in U_F(R(T_0))\cap U_F(F_1)), \cdots ,  S_{n+1}\in U_F(F_n)\cap U_F(R(T_*)). \eqno (2.12)$$
We define inductively
$$T_{m+i}=P_{F_i}^{S_i}T_{m+i-1} ,i=1, \cdots, n.$$ By induction, it is Evident
$$
N(T_{m+i})=N(T_m),  \mbox{and} \, R(T_{m+i})=F_i, i=1, \cdots, n.\eqno (2.13) .$$
By Theorem $1.3$ we have that $T_{m+i}\sim T_{m+i-1}$ in the set $\Omega_i=\{T\in B^+(E,F): R(T)\oplus S_i=F \, \mbox{and} \, N(T)(= N(T_{m+i-1}))=N(T_m)\} \, \, i=1, \cdots , n .$ Obviously the following results for $T\in \Omega_i $ hold : $ R_1\in U_E(N(T_0))\cap U_E(N_1), \cdots, U_E(N_{m-1})\cap U_E(N(T)$ (note $N(T))=N(T_m)=N_m)$;  $ S_1 \in U_F(R(T_0)\cap U_F(F_1), \cdots, S_i\in U_F(F_{i-1})\cap U_F(R(T)).$ This infers  $ \Omega_i \subset \widetilde{T_0} $ for $ i=1, \cdots , n ,$ and so
 $$T_0\sim T_{m +i}\, \, \mbox{in} \, \, \widetilde{T_0} \, \, \mbox{for} \,\, i=1, \cdots, n. \eqno (2.14)$$
 So the proof of the theorem is turned to that of  $$  T_{m+n}\sim T_*  \,  \, in \, \, \widetilde{{T_0}}.\eqno (2.15)$$
 By the assumptions $(2.9)$ and $(2.12)$, We have
$$N(T_{m+n})\oplus R_{m+1}=T_*\oplus R_{m+1}=E$$  $\mbox{and}$ $$ R(T_{m+n})\oplus S_{n+1}= R(T_*)\oplus S_{n+1}=F \eqno (2.16)$$
( note $ N(T_{m+n})=N_m \, \mbox{and} \, R(T_{m+n})=F_n $ ).
 Then by Theorem $1.4$, $T_* \sim  T_* P_{R_{m+1}}^{N(T_m)} $ in the set $ \Omega =\{T\in B^+(E,F): E=N(T) \oplus R_{m+1} \, \mbox{and} \, R(T)=R(T_*)\} .$ Let $\omega (t)$ be such a path in $\Omega$ that $\omega (0)=T_* \, \mbox{and} \, \omega (1)= P_{R_{m+1}}^{N(T_m)} T_*.$ Clearly $$ R_1\in U_E(R(T_0))\cap U_E(N_1), \cdots , N_{m+1}\in  U_E(N_m)\cap N(\omega(t))\eqno (2.17)$$ and $R(\omega (t))$ satisfies $(2.12).$ So $ T_* \sim P_{R_{m+1}}^{N(T_m)} T_*$ in $\widetilde{T_0}.$
Let $\gamma (t)= P_{F_n}^ {S_{n+1}} \omega (t).$ Because $ S_{n+1}\oplus R(T_*)=F \mbox{and} \, R(\omega (t))=R(T_*)$, $$ \gamma (t)x=0\Leftrightarrow \omega (t)x \in S_{n+1} \Leftrightarrow x\in N(\omega (t));$$ $$R(\gamma (t))= P_{F_n}^ {S_{n+1}}R(\omega (t))=P_{F_n}^ {S_{n+1}} (F_n\oplus S_{n+1})
=P_{F_n}^ {S_{n+1}}F=F_n.$$ Then since $ S_1\in U_F(R(T_0))\cap U_F(F_1)), \cdots , S_n\in U_F(F_{n-1})\cap U_F(\gamma (t))(=U_F(F_n))$  $\mbox{and} \, \mbox{and} \,$ the equality $N(T_{m+n})\oplus R_{m+1}=T_*\oplus R_{m+1}=E$ in $(2.16)$ we infer $$ T_*\sim P_{F_n}^ {S_{n+1}} T_* P_{R_{m+1}}^{N(T_m)} \, \,  \mbox{in} \, \, \widetilde{T_0} .$$  Let
$$T_{m+n}^+y=\left\{\begin{array}{rcl} (T_{m+n}\mid _ {R_{m+1} } )^{-1}y, & & {y\in F_n}\\
0, & & {y\in S_{n+1}}.
\end{array}\right.$$  Then $$P_{F_n}^ {S_{n+1}} T_* P_{R_{m+1}}^{N(T_m)}=P_{F_n}^ {S_{n+1}} T_*T_{m+n}^+T_{m+n}$$
( note $N(T_m)(=N(T_{m+n}))\oplus R_{m+1}=E \, \mbox{and} \, R(T_{m+n})(=F_n)\oplus S_{n+1}=F$ ).
We claim $$P_{F_n}^ {S_{n+1}} T_*T_{m+n}^+\mid_{F_n} \in B^\times(F_n),$$
Indeed, $$P_{Fn}^ {S_{n+1}} T_*T_{m+n}^+\mid_{F_n}y=0 \, \mbox{for} \, y\in F_n \Leftrightarrow T_*T_{m+n}^+y=0 \, \mbox{by} (2.16)\Leftrightarrow T_{m+n}^+y=0 \Leftrightarrow y=0 ; $$ while by $(2.16)$  we have for any $y\in F_n$ $$y=P_{Fn}^ {S_{n+1}}y =P_{Fn}^ {S_{n+1}} T_*r \, \mbox{for some} \, r\in R_{m+1} ,$$ and moreover , $$ y=
 P_{Fn}^ {S_{n+1}} T_*T_{m+n}^+y_0 \, \mbox{for some} \, y_0\in F_n.$$
It is well known that $P_{F_n}^ {S_{n+1}} T_*T_{m+n}^+\mid_{F_n} \sim  \, \mbox{either} \, I_{F_n} \, \mbox{or} \, -I_{F_n} \, \mbox{in} \, B^\times (F_n),$  Say that $\omega (t)$ lies in $B^\times (F_n) \, \mbox{and}$ satisfies $\omega (0)=  P_{F_n}^ {S_{n+1}} T_*T_{m+n}^+\mid_{F_n} \, \mbox{and} \, \omega (1)= \,  \mbox{either} \, \, I_{F_n} \,  \mbox{or} \, -I_{F_n}.$   Let $\gamma (t)=\omega (t)T_{m+n}.$ Then $N(\gamma (t))=N(T_{n+m})=N(T_m) , \, \, \mbox{and} \,  \, R(\gamma(t))=F_n , $ so that $\gamma (t)\subset \widetilde{T_0},$
$\gamma (0)=  P_{F_n}^ {S_{n+1}} T_*T_{m+n}^+T_{m+n}=P_{F_n}^ {S_{n+1}} T_* P_{R_{m+1}}^{N(T_m)}, $ and $$ \gamma (1)= T_{m+n} \, \, \mbox{or} \, \, -T_{m+n}.$$ Hence $$P_{F_n}^ {S_{n+1}} T_* P_{R_{m+1}}^{N(T_m)}\sim T_{m+n} \, \, \mbox{or} \, \, -T_{m+n} \, \, \mbox{in} \, \, \widetilde{T_0}.$$  By Theorem $2,1$  we conclude that in either way , $ T_*\sim T_{m+n}$ in $\widetilde{T_0}.$ For the end of the proof of the theorem we have to show that in the case of $\mbox{dim}N(T_0)>0 \, \mbox{and codim}R(T_0)=0,$ $T_m\sim T_*$ in $\widetilde{T_0}.$ By $(2.10)$ and $(2.9)$,
 $$R(T_m)=R(T_0)=F , N(T_m)=N_m \, \mbox{and} \, \,R_{m+1}\oplus N(T_m)=R(T_{m+1}) \oplus N(T*)=E.$$ Then by Theorem $1.4$,
 $$ T_*=T_*P^{N(T_*)}_{R_{m+1}}\sim T_*P^{N(T_m)}_{R_{m+1}} \quad \mbox{in} \, \,  \widetilde{T_0},$$ say $\omega (t)$ is such a path with $N(\omega (t))\oplus R_{m+1}=E \mbox{and} \, R(\omega (t)))=R(T_0)$ that $\omega (0))= T_* \, \mbox{and} \, \omega (1)= T_*P^{N(T_m)}_{R_{m+1}}$ , clearly $\omega (t)\subset \widetilde{T_0}.$ Let $T_m^+=(T_m\mid R_{m+1})^{-1}$ then $T_*\sim T_*T_m^+T_m \, \mbox{in} \, \widetilde{T_0}.$ It is evident that $T_*T_m^+ \in B^\times (F),$ so that $T_*\sim T_m \, \mbox{or} \, -T_m $ in $\widetilde{T_0}.$  By Theorem $2.1$ we conclude that in either way $$T_*\sim T_m \sim T_0 \, \, \mbox{in} \, \, \widetilde{T_0}$$ . \quad $\Box$

{\bf Remark:}

(i) if $T_0 \in F_k$ and $k<\min\{ \dim E, \dim F\},$ then $\widetilde{T_0}=F_k;$

(ii) if $T_0 \in \Phi_{m,n}, \, m, n<\infty$ and either $m>0$ or $n>0$, then $\widetilde{T_0}=\Phi_{m,n}.$

These are immediate from Theorems $2.2$ and $2,3.$

\vskip 0.2cm\begin{center}{\bf 3\quad Some Applications to Topology and Analysis in operators}
\end{center}\vskip 0.2cm

In $2007$ the following result was presented in [Ma 5]: both $F_k$ $(k<\infty ), \Phi_{m,n} (m,n<\infty)$ are smooth submanifolds in $B(E,F)$ with the tangent space $\{ T\in B(E,F): TN(X)\subset R(X)\} $ at any $X$ in them. Now it can be consolidated as follows,

{\bf Theorem 3.1}\quad {\it Both $F_k (k<\infty ),$ $\Phi_{m,n}$ with either $m>0$ or $n>0$ are smooth and path connected submanifolds in $B(E,F)$ with the tangent space $\{T\in B(E,F): TN(X)\subset R(X)\}$ at any $X$ in them .}\\
This is immediate from Theorems $2.2$ and $2.3.$

{\bf Theorem 3.2} \quad {\it $B(\mathbf{R}^m,\mathbf{R}^n)=\bigcup^{\min\{n,m\}}\limits_{k=0}F_k \,
, \, \, F_k (k<\min \{m,n\})$ is a smooth and path connected submanifold in
$B(\mathbf{R}^m,\mathbf{R}^n)$, and especially, $\dim F_k=(m+n-k)k$
for $k=0,1,\cdots,\min\{n,m\}.$}

{\bf Proof } \quad In the case of $k<\min \{m,n\}$, $T\in F_k$ implies either $ \dim N(T)>0$ or $\mathrm{codim} R(T)>0,$ so that it follows that $F_k$ is a smooth and path connected hyper surface in $B(\mathbf {R}^m, \mathbf{ R}^n)$ from theorem $3.1$.The Connectedness of $F_k$ allows us to take such a
special point $A\in F_k$ that $\mbox{dim} T_AF_k$ can be calculated. Let $ T=\{ t_{i,j}\}_{i,j=1}^{m,n}$  denote a matrix in $B(\mathbf {R}^m, \mathbf{ R}^n),$ and take $ A =\mathbf {I}(m \times m)$ where  $$ \mathbf {I}(m \times m) =  \left\{ \{t{i,j}\}_{i,j=1}^{m,n} : t_{i,j} =0 \, \, \mbox {except } \, \, t_{i,i}=1 , \, \, 1\leq i \leq k \right\} .$$ Next go to determine $N(A) , \, \, R(A)$ and $T_AF_k.$ Obviously $$N(A)= \{ (x_1, x_2, \, . \, . \, . ,x_m)\in R^m: x_1=x_2= .\, . \, . =x_k=0 \}$$ and $$ R(A)=\{ (y_1, y_2, .\, . \, . y_n)\in R^n: y_{k+1}=y_{k+2}= \, . \, . \, .=y_n=0\} \eqno (3.1) .$$ Directly $$ T= \mathbf{I}_k(n)T + \mathbf {I}_k^+T\mathbf{I}(m) + \mathbf{I}_k^+(n)T\mathbf{I}_k^+(m), \eqno(3.2)$$ where $\mathbf {I}_k(n)=\mathbf{I}(n\times n),\, \mathbf{I}_k(m)=\mathbf{I}_k(m\times m),$  $$\mathbf {I}_k^+(n)=\left \{ \{t_{i,j}\}_{i,j=1}^n \in B(\mathbf{R}^n): t_{i,j}=0 \, \, \, \mbox{except} \, \, \, t_{i,i}=1 , \, k+1\leq i\leq n\right\}, $$ and  $$ \mathbf{I}_k^+(m)=\left\{ \{t_{i,j}\}_ {i,j=1}^m \in B(\mathbf{R}^m): t_{i,j}=0 \, \, \, \mbox{except} \, \, \ \, t_{i,i}=1, \, k+1\leq i \leq m \right\}. $$ Note \, $\mathbf {I}_k^+(n) +\mathbf {I}_k (n)=\mathbf {I}_n \, $(the identity on $B(\mathbf{R}^n ) )$ and   $\mathbf {I}_k^+(m) +\mathbf {I}_k (m)=\mathbf {I}_m \, $(the identity on $B(\mathbf{R}^m ) )$. By $(3.1)$ and $(3.2)$ one observes
$$T_A F_k=\left \{ \mathbf{I}_k(n) T + \mathbf{I}_k^+ (n) T \mathbf{I}_k(m): \forall T\in B(\mathbf{R}^m ,\mathbf{R}^n ) \right \}. \quad \eqno(3.3)  $$
This is immediate from $(3.1)$ and $(3.2).$ Let $B=\{b_{i,j}\}_{i,j=1}^{m,n}$ for any $B\in T_A F_k.$ Then we have $b_{i,j} =0 $ except
 $$ b_{i,j}=\left\{\begin{array}{rcl}  t_{i,j}, & & i=1, \cdots, k \, \, \mbox{and} \, \, j=1, \cdots, n \\
t_{i,j}, & & i=k+1, \cdots, m \, \, \mbox{and} \, \, j=i, \cdots, k.
\end{array}\right.$$
Therefore
\begin{eqnarray*}  \dim F_k  & = & n\times k + (m-k)k \\&=& (m+n -k)k \end{eqnarray*}
 $k=1, \cdots, \, \,\min  \{ m, n\}.$
The proof ends .\quad $\Box$

Of special interest is the dimensional formula of $F_k \, \, k=0,1, . . . , \min\{m.n\},$ which is a new result in algebraic geometry.
\newpage
\begin{center}{\bf References}
\end{center}
\vskip -0.1cm
\medskip
{\footnotesize
\def\REF#1{\par\hangindent\parindent\indent\llap{#1\enspace}\ignorespaces}

\REF{[Abr]}\ R. Abraham, J. E. Marsden, and T. Ratin, Manifolds,
tensor analysis and applications, 2nd ed., Applied Mathematical
Sciences 75, Springer, New York, 1988.

\REF{[An]}\ V. I. Arnol'd, Geometrical methods in the theory of
ordinary differential equations, 2nd ed., Grundlehren der
Mathematischen Wissenschaften 250, Springer, New York, 1988.

\REF{[Bo]}\ B. Booss,D.D.Bleecker, Topology and Analysis: the
Atyah-Singer Index Formula and Gauge-Theoretic physics, New York:
Springer-Verlag 1985.

\REF{[Caf]}\ V. Cafagra, Global invertibility and finite
solvability, pp. 1-30 in Nonlinear functional analysis (New York,
NJ, 1987), edited by P. S. Milojevic, Lecture Notes in Pure and
Appl. Math. 121, Dekker, New York, 1990.

\REF{[Ma1]}\ Jipu Ma, Dimensions of Subspaces in a Hilbert Space and
Index of Semi-Fredholm Operators, Sci, China, Ser. A 39:12(1986).

\REF{[Ma2]}\ Jipu Ma, Generalized Indices of Operators in B(H), Sci,
China, Ser. A , 40:12(1987).

\REF{[Ma3]}\ Jipu Ma, Complete Rank Theorem in Advanced Calculus and
Frobenius Theorem in Banach Space, arxiv:1407.5198v5[math.FA]23 Jan.
2015.

\REF{[Ma4]}\ Jipu Ma, Frobenius Theorem in Banach Space,  arXiv:
submit/1512157[math.FA] 19 Mar 2016.

\REF{[Ma5]}\ Jipu Ma, Three classes of smooth Banach manifold in B(E,F ),
Sci.China Ser. A, 50 (19) (2007) 1233-1239.

\textbf{Acknowledgments}\quad

I may have the honor to dedicate this paper to my teacher Professor Xu Men-ying . I will always be grateful for her supplying financial strain for
my parents with regular funds during my studies in Nanjing University. Without this help I would have been unable to continue my studies in the
 university.

I would like to thank National Natural Science Foundation of China for continuous support of  my work in the past time .

\end{document}